\documentclass[a4paper,11pt]{article}
\usepackage{amsmath, amssymb, amsthm}
\usepackage{graphicx} 

\begin{document}

\def\Sym{{\rm Sym}} \def\Aut{{\rm Aut}} \def\Out{{\rm Out}}
\def\Id{{\rm Id}} \def\Z{{\mathbb Z}} \def\R{{\mathbb R}}
\def\N{{\mathbb N}} \def\Pc{{\rm Pc}} \def\Ker{{\rm Ker}}

\title{\bf{Irreducible Coxeter groups}}
 
\author{\textsc{Luis Paris}}

\date{\today}

\maketitle

\begin{abstract}

\noindent
We prove that a non-spherical irreducible Coxeter group is (directly) indecomposable and that a 
non-spherical and non-affine Coxeter group is strongly indecomposable in the sense that all its finite 
index subgroups are (directly) indecomposable.

\bigskip\noindent
Let W be a Coxeter group. Write $W= W_{X_1} \times \dots \times W_{X_b} \times W_{Z_3}$, where 
$W_{X_1}, \dots, W_{X_b}$ are non-spherical irreducible Coxeter groups and $W_{Z_3}$ is a finite one. 
By a classical result, known as the Krull-Remak-Schmidt theorem, the group $W_{Z_3}$ has a 
decomposition $W_{Z_3}= H_1 \times \dots \times H_q$ as a direct product of indecomposable groups, 
which is unique up to a central automorphism and a permutation of the factors. Now, $W= W_{X_1} \times 
\dots \times W_{X_b} \times H_1 \times \dots \times H_q$ is a decomposition of $W$ as a direct product 
of indecomposable subgroups. We prove that such a decomposition is unique up to a central automorphism and 
a permutation of the factors.

\bigskip\noindent
Write $W=W_{X_1} \times \dots \times W_{X_a} \times W_{Z_2} \times W_{Z_3}$, where $W_{X_1}, \dots, 
W_{X_a}$ are non-spherical and non-affine Coxeter groups, $W_{Z_2}$ is an affine Coxeter group whose 
irreducible components are all infinite, and $W_{Z_3}$ is a finite Coxeter group. The group $W_{Z_2}$ 
contains a finite index subgroup $R$ isomorphic to $\Z^d$, where $d=|Z_2|-b+a$ and $b-a$ is the number 
of irreducible components of $W_{Z_2}$. Choose $d$ copies $R_1, \dots, R_d$ of $\Z$
such that $R= R_1 \times \dots \times R_d$. Then $G=W_{X_1} \times \dots \times W_{X_a} \times R_1 
\times \dots \times R_d$ is a virtual decomposition of $W$ as a direct product of strongly indecomposable subgroups. 
We prove that such a virtual decomposition is unique up to commensurability and a permutation of the 
factors.

\end{abstract}


\section{Introduction}

Let $S$ be a finite set. A {\it Coxeter matrix} over $S$ is a square matrix $M=(m_{s\,t})_{s,t\in 
S}$ indexed by the elements of $S$ and such that $m_{s\,s}=1$ for all $s \in S$ and 
$m_{s\,t}=m_{t\,s} \in \{2,3,4, \dots, +\infty\}$ for all $s,t \in S$, $s \neq t$. A Coxeter 
matrix $M=(m_{s\,t})_{s,t\in S}$ is usually represented by its {\it Coxeter graph}, $\Gamma$, 
which is defined by the following data. The set of vertices of $\Gamma$ is $S$, two vertices $s,t \in S$ are 
joined by an edge if $m_{s\,t} \ge 3$, and this edge is labeled by $m_{s\,t}$ if $m_{s\,t} \ge 4$.

\bigskip\noindent
Let $\Gamma$ be a Coxeter graph with set of vertices $S$. Define the {\it Coxeter system of type 
$\Gamma$} to be the pair $(W,S)$, where $W=W_\Gamma$ is the group generated by $S$ and subject to 
the relations
\begin{gather*}
s^2=1 \quad \text{for } s \in S\,,\\
(st)^{m_{s\,t}}=1 \quad \text{for } s,t \in S,\ s \neq t, \text{ and } m_{s\,t}<+\infty\,,
\end{gather*}
where $M=(m_{s\,t})_{s,t\in S}$ is the Coxeter matrix of $\Gamma$.  The group $W=W_\Gamma$ is 
simply called the {\it Coxeter group of type $\Gamma$}.

\bigskip\noindent
A Coxeter system $(W,S)$ is called {\it irreducible} if its defining Coxeter graph is connected. 
The number $n$ of elements of $S$ is called the {\it rank} of the Coxeter system. 

\bigskip\noindent
The term ``Coxeter'' in the above definition certainly comes from the famous Coxeter's theorem 
\cite{Coxet2}, \cite{Coxet1} which states that the finite (real) reflection groups are precisely the 
finite Coxeter groups. However, the Coxeter groups were first introduced by Tits in an unpublished 
manuscript \cite{Tits1} whose results appeared in the seminal Bourbaki's book \cite{Bourb1}. The 
Coxeter groups have been widely studied, they have
many attractive properties, and they form an important source of examples for 
group theorists. Basic references for them are \cite{Bourb1} and \cite{Humph1}.

\bigskip\noindent
The notion of irreducibility is of importance in the theory as all the Coxeter groups can be 
naturally decomposed as direct products of irreducible ones. Nevertheless, I do not 
know any work which studies the irreducibility itself. This is the object of the present paper. 

\bigskip\noindent
Our analysis needs to differentiate three classes of irreducible systems: the spherical ones which 
correspond to the finite Coxeter groups, the affine ones, and the remainder that we shall simply call 
non-spherical and non-affine Coxeter systems. These classes are defined as follows.

\bigskip\noindent
Let $\Gamma$ be a Coxeter graph and let $(W,S)$ be the Coxeter system of type $\Gamma$. Take an 
abstract set $\Pi=\{ \alpha_s; s \in S\}$ in one-to-one correspondence with $S$ and call the 
elements of $\Pi$ {\it simple roots}. Let $V= \oplus_{s\in S} \R \alpha_s$ be the real vector 
space having $\Pi$ as a basis. Define the {\it canonical form} to be the symmetric bilinear form 
$\langle\cdot,\cdot \rangle : V \times V \to \R$ determined by
\[
\langle \alpha_s, \alpha_t \rangle= \left\{
\begin{array}{ll}
-\cos (\pi/m_{s\,t}) \quad &\text{if }m_{s\,t} <+\infty\,,\\
-1\quad &\text{if }m_{s\,t}= +\infty\,,
\end{array}\right.
\]
for all $s,t \in S$. There is a faithful linear representation $W \to GL(V)$, called the {\it canonical 
representation}, which leaves invariant the canonical form, and which is defined by
$s(x)= x-2\langle x,\alpha_s \rangle \alpha_s$,
for $x \in V$ and $s \in S$.

\bigskip\noindent
A Coxeter graph $\Gamma$ is of {\it spherical type} if the canonical form of $\Gamma$ is 
positive definite. As pointed out before, it has been prove by Coxeter 
\cite{Coxet2}, \cite{Coxet1} that $\Gamma$ is of spherical 
type if and only if the associated Coxeter group $W_\Gamma$ is finite. The connected spherical type 
Coxeter graphs are precisely the Coxeter graphs $A_n$ ($n\ge 1$), $B_n$ ($n\ge 2$), $D_n$ ($n\ge 4$), 
$E_6$, $E_7$, $E_8$, $F_4$, $G_2$, $H_3$, $H_4$, $I_2(p)$ ($p\ge 5$ and $p \neq 6$) pictured in 
\cite{Bourb1}, p. 193.

\bigskip\noindent
A Coxeter graph $\Gamma$ is of {\it affine type} if the canonical form of $\Gamma$ is 
positive semidefinite but not positive definite. Any affine Coxeter group $W_\Gamma$ can be presented 
as a semidirect product $R \rtimes W_0$, where $R$ is a finitely generated free abelian group of rank $\ge 1$ and 
$W_0$ is a finite Coxeter group. If, moreover, $\Gamma$ is connected, then the rank of $R$ is precisely $n-1$, 
where $n=|S|$ is the rank of $(W,S)$. Conversely, if a Coxeter group $W_\Gamma$ 
has a finite index subgroup isomorphic to $\Z^d$ ($d\ge 1$), then $\Gamma$ is of affine type (see 
\cite{Harpe1}). The connected affine type Coxeter graphs are precisely the Coxeter graphs $\tilde A_n$ ($n \ge 
1$), $\tilde B_n$ ($n \ge 2$), $\tilde C_n$ ($n\ge 3$), $\tilde D_n$ ($n\ge 4$), $\tilde E_6$, $\tilde E_7$, 
$\tilde E_8$, $\tilde F_4$, $\tilde G_2$ pictured in \cite{Bourb1}, p. 199.

\bigskip\noindent
Let $G$ be a group. A subgroup $H$ of $G$ is called a {\it direct factor} of $G$ if there exists a 
subgroup $K$ of $G$ such that $G= H \times K$. If there are no proper non-trivial direct factors of 
$G$, then $G$ is said to be {\it indecomposable}. If $G$ is infinite and all the finite index subgroups 
of $G$ are indecomposable, then $G$ is said to be {\it strongly indecomposable}.

\bigskip\noindent
Clearly, every simple group is indecomposable. Other examples include the infinite cyclic group, the 
cyclic group of order $p^m$, where $p$ is a prime number, the free groups, and the Artin groups of 
spherical type \cite{Paris1}. Now, the infinite simple groups are also strongly indecomposable as well as 
$\Z$ and the torsion free hyperbolic groups. In this paper we prove that an irreducible non-spherical 
Coxeter group is indecomposable and that an irreducible non-spherical and non-affine Coxeter group is 
strongly indecomposable (see Theorem 4.1). (After the first version of the present paper, the first 
part of Theorem 4.1 has been extended to the infinitely generated Coxeter groups by Nuida 
\cite{Nuida1}, and a different proof of Theorem 4.1 has been proposed by de Cornulier and de la Harpe 
\cite{CorHar1}.) 

\bigskip\noindent
Not all the irreducible spherical type Coxeter groups are indecomposable.
The reader certainly 
knows the example of the dihedral group of order $8q+4$ which is isomorphic to the direct product 
of $C_2= \{ \pm 1 \}$ with the dihedral group of order $4q+2$. However, there are few exceptions. 
In particular, one of the factors must be the center of the group. The decompositions as direct 
products of the finite irreducible Coxeter groups are known to experts, 
but I do not know any reference where they are listed. So, we include the list in the last section
of the paper.

\bigskip\noindent
An irreducible affine type Coxeter group $W_\Gamma$ cannot be strongly indecomposable since it contains 
a finite index subgroup isomorphic to $\Z^{n-1}$, where $n$ is the rank of $W_\Gamma$, except if $n=2$. 
If $n=2$, then $W_\Gamma$ is the Coxeter group of type $\tilde A_1$ which is isomorphic to $C_2 \ast 
C_2$, where $C_2=\{ \pm 1\}$.

\bigskip\noindent
A {\it Remak decomposition} of a group $G$ is a decomposition of $G$ as a direct product 
of finitely many non-trivial indecomposable subgroups. 
An automorphism $\varphi: G \to G$ is said to be {\it central} if it is the identity modulo the 
center of $G$, that is, if $g^{-1} \varphi(g) \in Z(G)$ for all $g \in G$.
We say that a group $G$ satisfies the {\it 
maximal} (resp. {\it miniùmal}) {\it condition on normal subgroups} if each non-empty family of normal 
subgroups contains at least one maximal (resp. minimal) element for the inclusion. Obviously, finite 
groups satisfy both conditions. On the other hand, the Coxeter groups do not satisfy these conditions 
in general. For example, the Coxeter group $C_2 \ast C_2 \ast C_2$ does not satisfy neither the maximal 
condition, nor the minimal condition on normal subgroups. The proof of this last fact is left to the 
reader (but can be also found in \cite{CorHar1}).

\bigskip\noindent
A classical result, known as the Krull-Remak-Schmidt theorem, states that, if a group $G$ satisfies 
either the maximal or the minimal condition on normal subgroups, then it has a Remak decomposition. 
Furthermore, if it satisfies both, the maximal and the minimal conditions on normal subgroups, then 
this Remak decomposition is unique up to a central automorphism and a permutation of the factors. Of 
course, there are groups which have no Remak decompositions, and there are groups which have two 
different Remak decompositions whose factors are not isomorphic up to a permutation. We refer to 
\cite{CorHar1} for a detailed account on this question. Our aim here is to show that a Coxeter group 
has a Remak decomposition which is unique up to a central automorphism (see Theorem 5.2). This Remak 
decomposition is quite natural and can be described as follows.

\bigskip\noindent
Let $\Gamma$ be a Coxeter graph and let $(W,S)$ be the Coxeter system of type $\Gamma$. For $X 
\subset S$, we denote by $\Gamma_X$ the full subgraph of $\Gamma$ generated by $X$ and by $W_X$ 
the subgroup of $W$ generated by $X$. By \cite{Bourb1}, Chap. 4, \S\ 8, $(W_X,X)$ is the Coxeter system of type 
$\Gamma_X$. A subgroup of the form $W_X$ is called a {\it standard parabolic subgroup}, and a 
subgroup conjugated to a standard parabolic subgroup is simply called a {\it parabolic subgroup}. 

\bigskip\noindent
Let $\Gamma_1, 
\dots, \Gamma_l$ be the connected components of $\Gamma$ and, for $1 \le i\le l$, let $X_i$ be the set 
of vertices of $\Gamma_i$. Then $W=W_{X_1} \times \dots \times W_{X_l}$ and each $(W_{X_i}, X_i)$ 
is an irreducible 
Coxeter system. The above equality is called the {\it standard decomposition of $W$ into irreducible 
components}, and each $W_{X_i}$ is called an {\it irreducible component} of $W$.
Up to a permutation of the $\Gamma_i$'s, we can assume that, for some $0 \le 
a\le b\le l$, the Coxeter graphs $\Gamma_1, \dots, \Gamma_a$ are non-spherical and non-affine,
the Coxeter graphs $\Gamma_{a+1}, \dots, \Gamma_b$ are of affine type, and the Coxeter 
graphs $\Gamma_{b+1}, \dots, \Gamma_l$ are of spherical type. Set $Z_1= X_1 \cup \dots \cup X_a$, 
$Z_2= X_{a+1} \cup \dots \cup X_b$, and $Z_3= X_{b+1} \cup \dots \cup X_l$. Then $W=W_{Z_1} \times 
W_{Z_2} \times W_{Z_3}$, $W_{Z_1}$ is a Coxeter group whose irreducible components are all 
non-spherical and non-affine,
$W_{Z_2}$ is an affine Coxeter group, and $W_{Z_3}$ is a finite Coxeter group. The subgroup 
$W_{Z_1}$ is called the {\it non-spherical and non-affine part} of $W$, 
$W_{Z_2}$ is called the {\it affine part}, 
and $W_{Z_3}$ is called the {\it finite part} of $W$.

\bigskip\noindent
By the above mentioned Krull-Remak-Schmidt theorem, the finite part has a Remak decomposition $W_{Z_3}= 
H_1 \times \dots \times H_q$ which is unique up to a central automorphism. Then
\[
W=W_{X_1} \times \dots \times W_{X_b} \times H_1 \times \dots \times H_q
\]
is a Remak decomposition of $W$.

\bigskip\noindent
Let $G$ be a group. Define a {\it virtual strong Remak decomposition} of $G$ to be a finite index 
subgroup $H$ provided with a decomposition $H= H_1 \times \dots \times H_m$ as a direct product 
of finitely many (infinite) strongly indecomposable subgroups. Recall that two groups $G_1$ and $G_2$ are 
{\it commensurable} if there is a finite index subgroup of $G_1$ isomorphic to a finite 
index subgroup of $G_2$. 
Two virtual strong Remak decompositions $H= H_1 \times \dots \times H_p$ and $H'= H_1' \times 
\dots \times H_q'$ of $G$
are called {\it equivalent} if
$p=q$ and, up to a permutation of the $H_i$'s, $H_i$ and $H_i'$ are 
commensurable for all $1 \le i\le p$.

\bigskip\noindent
Obviously, a strong indecomposable group as well as a finitely generated abelian group has a unique 
virtual strong Remak decomposition up to equivalence. I do not know much more about virtual strong 
Remak decompositions, but I believe that this question should be of interest in the study of groups up 
to commensurability. We prove in this paper (Theorem 6.1) that a Coxeter group has a unique virtual 
strong Remak decomposition up to equivalence. Here again, the decomposition is quite natural and can be 
described as follows.

\bigskip\noindent
Let $\Gamma$ be a Coxeter graph and let $(W,S)$ be the Coxeter system of type $\Gamma$. Write 
$W=W_{X_1} \times \dots \times W_{X_a} \times W_{Z_2} \times W_{Z_3}$, where $W_{X_1}, \dots, W_{X_a}$ 
are the non-spherical and non-affine irreducible components of $W$, $W_{Z_2}$ is the affine part, and 
$W_{Z_3}$ is the finite part. Let $d= |Z_2| -b+a$, where $b-a$ is the number of affine irreducible 
components of $(W,S)$. By \cite{Bourb1}, Chap. 6, \S\ 2, $W_{Z_2}$ contains a finite index subgroup $R$ 
isomorphic to $\Z^d$. Let $R_1, \dots, R_d$ be $d$ copies of $\Z$ such that $R=R_1 \times \dots \times 
R_d$. Then
\[
G= W_{X_1} \times \dots \times W_{X_a} \times R_1 \times \dots \times R_d
\]
is a virtual strong Remak decomposition of $W$. The number $d$ will be called the {\it affine 
dimension} of $W$.

\bigskip\noindent
An element $w \in W$ is called {\it essential} if it does not lie in any proper parabolic 
subgroup.
Part of the proofs of the paper
are based on Krammer's study on essential elements (see \cite{Kramm1}). I am 
sure that other proofs with other (more or less simple) techniques can be easily found by 
experts. However, another goal of the present paper is to present this piece of Krammer's Ph. D. Thesis 
which deserves to be known. In particular, it is shown in \cite{Kramm1} that, in a non-spherical and 
non-affine Coxeter group $W$,
\begin{itemize}
\item 
every essential element is of infinite order,
\item 
an element $w$ is essential if and only if any non-trivial power $w^m$ is essential,
\item 
the centralizer in $W$ of an essential element $w$ contains $\langle w \rangle = \{ w^m; m 
\in \Z \}$ as a finite index subgroup.
\end{itemize}
These Krammer's results are explained in Section 2.

\bigskip\noindent
Choose a linear ordering $S=\{ s_1, s_2, \dots, s_n\}$ and define the {\it Coxeter element} of 
$W$ (with respect to this linear ordering) to be $c=s_n \dots s_2s_1$. In order to apply the 
above mentioned Krammer's results,
we need at some point to choose an essential element. In particular, we first have to show that 
such an element exists. A natural approach to this question is to prove that the Coxeter 
elements are essential. Curiously, this result is unknown (for instance, it is a question in \cite{McMul1}). 
We give a simple proof of this fact in Section 3.


\section{Essential elements}

Let $\Gamma$ be a Coxeter graph and let $(W,S)$ be the Coxeter system of type $\Gamma$. Recall from the 
introduction the set $\Pi= \{ \alpha_s; s \in S\}$ of simple roots, the real vector space $V= 
\oplus_{s \in S} \R \alpha_s$, the canonical form $\langle\cdot, \cdot \rangle : V \times V \to \R$, and the 
canonical representation $W \to GL(V)$.
The set $\Phi= \{ w\alpha_s; w \in W \text{ and } s \in S\}$ is called the {\it root system} of 
$(W,S)$. The elements of $\Phi^+= \{ \beta= \sum_{s \in S} \lambda_s \alpha_s \in \Phi; \lambda_s 
\ge 0 \text{ for all } s \in S\}$ are the {\it positive roots}, and the elements of $\Phi^-= -
\Phi^+$ are the {\it negative roots}. For $w \in W$, we set
$\Phi_w= \{ \beta \in \Phi^+; w \beta \in \Phi^- \}$.
The following proposition collects some classical results (see for example Chapter 5 of \cite{Humph1}).

\bigskip\noindent
{\bf Proposition 2.1.} {\it 
\begin{enumerate}
\item We have the disjoint union $\Phi= \Phi^+ \sqcup \Phi^-$.
\item Let $\lg: W \to \N$ denote the word length with respect to $S$. Let $w \in W$. Then $\lg(w)= 
|\Phi_w|$. In particular, $\Phi_w$ is finite.
\item Let $w \in W$ and $s \in S$. Set $\beta = w\alpha_s$ and $r_\beta= wsw^{-1}$. Then 
$r_\beta$ acts on $V$ by $r_\beta(x)= x -2 \langle x,\beta \rangle \beta$, for all $x \in V$.
\end{enumerate}}

\bigskip\noindent
Let $u,v \in W$ and $\alpha \in \Phi$. We say that $\alpha$ {\it separates} $u$ and $v$ if there 
is some $\varepsilon \in \{ \pm 1\}$ such that $u \alpha \in \Phi^\varepsilon$ and $v \alpha \in 
\Phi^{- \varepsilon}$. Let $w \in W$ and $\alpha \in \Phi$. We say that $\alpha$ is {\it 
$w$-periodic} if there is some $m \ge 1$ such that $w^m \alpha = \alpha$.

\bigskip\noindent
{\bf Proposition 2.2.} {\it Let $w \in W$ and $\alpha \in \Phi$. Then exactly one of the 
following holds.
\begin{enumerate}
\item $\alpha$ is $w$-periodic.
\item $\alpha$ is not $w$-periodic, and the set $\{ m \in \Z; \alpha \text{ separates } w^m \text{ 
and } w^{m+1} \}$ is finite and even.
\item $\alpha$ is not $w$-periodic, and the set $\{ m \in \Z; \alpha \text{ separates } w^m \text{ 
and } w^{m+1} \}$ is finite and odd.
\end{enumerate}}

\bigskip\noindent
We say that $\alpha$ is {\it $w$-even} is Case 2, and {\it $w$-odd} in Case 3.

\bigskip\noindent
{\bf Proof.} Assume that $\alpha$ is not $w$-periodic. Set $N(\alpha)= \{ m \in \Z; \alpha \text{ 
separates } w^m \text{ and } w^{m+1}\}$. Clearly, we only need to show that $N(\alpha)$ is 
finite.

\bigskip\noindent
Observe that, if $m \in N(\alpha)$, then $w^m \in \Phi_w \cup -\Phi_w$. On the other hand, if 
$w^{m_1} \alpha= w^{m_2} \alpha$, then $w^{m_1-m_2} \alpha = \alpha$, thus $m_1=m_2$, since 
$\alpha$ is not $w$-periodic. The set $\Phi_w \cup -\Phi_w$ is finite (see Proposition 2.1), thus 
$N(\alpha)$ is finite, too.
\qed

\bigskip\noindent
{\bf Proposition 2.3.} {\it Let $\alpha \in \Phi$, $w \in W$, and $p \in \N$, $ p\ge 1$.
\begin{itemize}
\item $\alpha$ is $w$-periodic if and only if $\alpha$ is $w^p$-periodic.
\item $\alpha$ is $w$-even if and only if $\alpha$ is $w^p$-even.
\item $\alpha$ is $w$-odd if and only if $\alpha$ is $w^p$-odd.
\end{itemize}}

\bigskip\noindent
{\bf Proof.} The equivalence
\[
\alpha \text{ is $w$-periodic} \quad \Leftrightarrow \quad \alpha \text{ is $w^p$-periodic}
\]
is obvious. So, we can assume that $\alpha$ is not $w$-periodic and we show the other 
equivalences.

\bigskip\noindent
Observe that, if $\alpha$ is $w$-even, then there are $m_0 >0$ and $\varepsilon \in \{ 
\pm 1\}$ such that $w^m \alpha, w^{-m} \alpha \in \Phi^\varepsilon$ for all $m \ge m_0$. 
On the other hand, if $\alpha$ is $w$-odd, then there are $m_0>0$ and $\varepsilon \in \{ \pm 
1\}$ such that $w^m \alpha \in \Phi^\varepsilon$ and $w^{-m} \alpha \in \Phi^{-\varepsilon}$ for 
all $m \ge m_0$. It is easily checked that these criteria for a non-periodic root to be even or 
odd imply the last two equivalences.
\qed

\bigskip\noindent
Recall that an {\it essential element} is an element $w \in W$ which does not lie in any proper
parabolic subgroup of $W$.
Now, we state the first Krammer's result in which we are interested (see \cite{Kramm1}, Corollary 
5.8.7).

\bigskip\noindent
{\bf Theorem 2.4 (Krammer \cite{Kramm1}).} {\it Assume that $\Gamma$ is connected and non-spherical.
Let $w \in W$. Then $w$ is essential if and only if $W$ is generated by the set $\{ 
r_\beta; \beta \in \Phi^+ \text{ and } \beta \text{ $w$-odd}\}$.}

\bigskip\noindent
{\bf Remark.} It will be shown in Section 3 that the Coxeter elements are essential, but the proof
will not use Theorem 2.4. Actually, I do not know how to use Theorem 2.4 in a simple way for this purpose.

\bigskip\noindent
A direct consequence of Proposition 2.3 and Theorem 2.4 is the following.

\bigskip\noindent
{\bf Corollary 2.5.} {\it Assume that $\Gamma$ is connected and non-spherical. Let $w \in W$ and 
$p \in \N$, $p \ge 1$. Then $w$ is essential if and only if $w^p$ is essential.}

\bigskip\noindent
Now, the following 
result is one of the main tools in the present paper. It 
can be found in \cite{Kramm1}, Corollary 6.3.10.

\bigskip\noindent
{\bf Theorem 2.6 (Krammer \cite{Kramm1}).} {\it Assume that 
$\Gamma$ is connected, non-spherical, and non-affine. 
Let $w \in W$ be an essential element. Then $\langle w \rangle = \{w^m; m \in \Z\}$ is a 
finite index subgroup of the centralizer of $w$ in $W$.}


\section{Coxeter elements}

{\bf Theorem 3.1.} {\it Let $\Gamma$ be a Coxeter graph and let $(W,S)$ be the Coxeter system of 
type $\Gamma$. Then any Coxeter element $c=s_n \dots s_2s_1$ of $W$ is essential.}

\bigskip\noindent
The following lemma is a preliminary result to the proof of Theorem 3.1.

\bigskip\noindent
{\bf Lemma 3.2.} {\it Let $\Gamma$ be a Coxeter graph with set of vertices $S$. Then there exist a 
Coxeter graph $\tilde \Gamma$ with set of vertices $\tilde S$ and an embedding $S \hookrightarrow 
\tilde S$ such that $\Gamma= \tilde \Gamma_S$ and the canonical form $\langle\cdot, 
\cdot\rangle_{\tilde \Gamma}$ of $\tilde \Gamma$ is non-degenerate.}

\bigskip\noindent
{\bf Proof.} We denote by $B_\Gamma = (b_{s\,t})_{s,t\in S}$ the Gram matrix of the canonical form 
$\langle\cdot,\cdot \rangle_\Gamma$ in the basis $\Pi= \{\alpha_s; s \in S\}$. So, $b_{s\,t}= -\cos( \pi/ 
m_{s\,t})$ if $m_{s\,t}< +\infty$, and $b_{s\,t}=-1$ if $m_{s\,t}=+\infty$. For $X \subset S$, we 
shall simply write $B_X$ for $B_{\Gamma_X}$.

\bigskip\noindent
We can and do choose a non-empty subset $X_0 \subset S$ such that $\det B_{X_0} \neq 0$, and 
such that $\det B_X=0$ whenever $|X|>|X_0|$ for $X \subset S$. Set $X_1=S \setminus X_0$, take an 
abstract set $\tilde X_1= \{ \tilde s; s \in X_1\}$ in one-to-one correspondence with $X_1$, and 
set $\tilde S= X_0 \sqcup X_1 \sqcup \tilde X_1$. Define $\tilde \Gamma$ adding to $\Gamma$ an 
edge labeled by $+\infty$ between $s$ and $\tilde s$ for all $s \in X_1$. We emphasize that 
there is no other edge in $\tilde \Gamma$ apart the ones of $\Gamma$, and the edges labeled by 
$+\infty$ between the elements of $X_1$ and they corresponding elements of $\tilde X_1$. Now, the 
proof of the following equality is left to the reader.
\[
\det B_{\tilde \Gamma} = (-1)^{|X_1|} \det B_{X_0} \neq 0\,.
\]
\qed

\bigskip\noindent
{\bf Proof of Theorem 3.1.} We suppose that $c=s_n \dots s_2s_1$ is not essential. So, there 
exist $X \subset S$, $X \neq S$, and $u \in W$, such that $c \in u W_X u^{-1}$.

\bigskip\noindent
Our proof uses the observation that, if $w \in W$ is not essential, then there exists a non-zero 
vector $x \in V \setminus \{0\}$ such that $w(x)=x$. The converse is false, especially when the 
canonical form $\langle\cdot, \cdot\rangle_\Gamma$ is degenerate. Indeed, all the 
vectors of the radical of $\langle\cdot,\cdot \rangle_\Gamma$ are fixed by all the elements of $W$. In 
order to palliate this difficulty, we use the simple trick to embed $W$ into a larger Coxeter 
group with a non-degenerate canonical form.

\bigskip\noindent
So, by Lemma 3.2, there exist a Coxeter graph $\tilde \Gamma$ with set of vertices $\tilde S$ and 
an embedding $S \hookrightarrow \tilde S$ such that $\tilde \Gamma_S = \Gamma$ and the canonical 
form $\langle\cdot,\cdot \rangle= \langle\cdot,\cdot \rangle_{\tilde \Gamma}$ of $\tilde \Gamma$ 
is non-degenerate. 
Let $\{s_{n+1}, \dots, s_m\}=\tilde S \setminus S$, let $\tilde W$ be the Coxeter group of type 
$\tilde \Gamma$, and let $\tilde V= \oplus_{i=1}^m \R \alpha_{s_i}$.

\bigskip\noindent
We write $\beta_i = u^{-1} (\alpha_{s_i})$ and $r_i= u^{-1} s_i u = r_{\beta_i}$ for $1 \le i\le 
n$, and we write $\beta_i = \alpha_{s_i}$ for $n+1 \le i\le m$. Note that $\{\beta_1, \dots, 
\beta_n, \beta_{n+1}, \dots, \beta_m \}$ is a basis of $\tilde V$. 
Let $d= s_m \dots s_{n+1} r_n \dots r_2r_1$. Observe that $d= s_m \dots s_{n+1} (u^{-1} c u)$, thus
$d \in \tilde W_{X\sqcup(\tilde S \setminus S)}$. Moreover, $\tilde W_{X \sqcup (\tilde S
\setminus S)}$ is a proper (standard) parabolic subgroup, hence
there exists $x \in \tilde V \setminus \{0\}$ such that $d(x)=x$. It is easily 
checked that there exist numbers $\lambda_{i\,j}= \lambda_{i\,j} (x) \in \R$, $1 \le i<j\le m$, 
such that
\[
d(x)=x-\sum_{j=1}^m \left( 2 \langle x, \beta_j \rangle + \sum_{i=1}^{j-1} \lambda_{i\,j} \langle 
x, \beta_i \rangle \right) \beta_j \,.
\]
Now, the equality $d(x)=x$ implies
\[
\sum_{j=1}^m \left( 2 \langle x, \beta_j \rangle + \sum_{i=1}^{j-1} \lambda_{i\,j} \langle x, 
\beta_i \rangle \right) \beta_j =0\,,
\]
thus, since $\{ \beta_1, \dots, \beta_m \}$ is a basis of $\tilde V$,
\[
\langle x, \beta_1 \rangle = \dots = \langle x, \beta_m \rangle = 0\,.
\]
This contradicts the fact that $x \neq 0$, $\{\beta_1, \dots, \beta_m \}$ is a basis, and 
$\langle\cdot,\cdot\rangle= \langle\cdot,\cdot\rangle_{\tilde \Gamma}$ is non-degenerate.
\qed


\section{ Indecomposability}

Let $\Gamma$ be a Coxeter graph and let $(W,S)$ be the Coxeter system of type $\Gamma$. The goal of the 
present section is to prove the following.

\bigskip\noindent
{\bf Theorem 4.1.}
{\it \begin{enumerate}
\item
Assume that $\Gamma$ is connected and non-spherical. Then $W$ is indecomposable.
\item
Assume that $\Gamma$ is connected, non-spherical, and non-affine. Then $W$ is 
strongly indecomposable.
\end{enumerate}}

\bigskip\noindent
We start collecting some well-known results on Coxeter groups. 
For $X \subset S$, we set $\Pi_X= \{ \alpha_s; s \in X\}$ and $G_X= \{w \in W; w\Pi_X=\Pi_X\}$. The set 
$T= \{ wsw^{-1}; s \in S \text{ and } w \in W\}$ is called the set of {\it reflections} of $W$. Any 
subgroup of $W$ generated by reflections is called a {\it reflection subgroup} of $W$.

\bigskip\noindent
For a group $G$ and a subgroup $H$ of $G$ we denote by $N_G(H)$ the normalizer of $H$ in $G$ and by 
$Z_G(H)$ the centralizer of $H$ in $G$. The center of $G$ is denoted by $Z(G)$.

\bigskip\noindent
{\bf Proposition 4.2.} {\it
\begin{enumerate}
\item{\rm (Bourbaki \cite{Bourb1}, Chap. 5, \S\ 4, Exercice 2).} Any finite subgroup of $W$ is 
contained in a finite parabolic subgroup.
\item
{\rm (Solomon \cite{Solom1}).} The intersection of a family of parabolic subgroups is a parabolic 
subgroup.
\item
{\rm (Deodhar \cite{Deodh1}).} If $X \subset S$, then $N_W(W_X)= W_X \rtimes G_X$.
\item
{\rm (Deodhar \cite{Deodh1}).} 
Let $T=\{w s w^{-1}; s \in S \text{ and } w \in W\}$ be the set of reflections. Then
$W$ is finite if and only if $T$ is finite.
\item
{\rm (Dyer \cite{Dyer1}).} If $H$ is a reflection subgroup of $W$, then there exists a subset $Y 
\subset H \cap T$ such that $(H,Y)$ is a Coxeter system.
\end{enumerate}}

\bigskip\noindent
Let $(W,S)$ be a non-spherical and non-affine Coxeter system, let $G$ be a 
finite index subgroup of $W$, and let $A,B$ be two subgroups of $G$ such that $G=A \times B$. The proof 
(of the second part) of Theorem 4.1 is divided into two steps. In the first one we show that either $A$ 
or $B$ is finite. The main tool in this step is Krammer's result stated in Theorem 2.6. 
The second step consists on showing that $A$ is trivial if it is finite. This second step is 
based on the following result.

\bigskip\noindent
{\bf Proposition 4.3.} {\it Assume that $\Gamma$ is connected and non-spherical. Let $H$ be a 
non-trivial finite subgroup of $W$. Then $N_W(H)$ has infinite index in $W$.}

\bigskip\noindent
{\bf Proof.}

\bigskip\noindent
{\bf Assertion 1.} {\it Let $X \subset S$, $\emptyset \neq X \neq S$. Then $W_X$ is not normal in $W$.}

\bigskip\noindent 
Let $V_X$ denote the linear subspace of $V$ spanned by $\Pi_X$. By Proposition 4.2, we have $w(V_X)=V_X$ for 
all $w \in N_W(W_X)$. Choose $s \in X$ and $t \in S \setminus X$ such that $m_{s\,t} \ge 3$. Then 
$\alpha_s \in V_X$ and $t(\alpha_s)= \alpha_s - 2\langle \alpha_s, \alpha_t \rangle \alpha_t \not\in V_X$, 
thus $t \not\in N_W(W_X)$.

\bigskip\noindent
{\bf Assertion 2.} {\it Let $H$ be a non-trivial finite subgroup of $W$. Then $H$ is not normal.}

\bigskip\noindent
Let $\Pc (H)$ denote the intersection of all the parabolic subgroups that contain $H$. By Proposition 
4.2, $\Pc (H)$ is a finite parabolic subgroup. Moreover, $\{1\} \neq \Pc (H)$ since it contains $H$, 
and $\Pc (H) \neq W$ since $\Pc (H)$ is finite while $W$ is infinite. By Assertion 1, if follows that 
$N_W(\Pc (H)) \neq W$. Now, observe that $N_W(H) \subset N_W( \Pc (H))$, thus $N_W(H) \neq W$.

\bigskip\noindent
{\bf Assertion 3.} {\it Recall that $T= \{wsw^{-1}; s \in S \text{ and } w \in W\}$. Let $X \subset S$, 
let $T_X= \{wxw^{-1}; x \in X \text{ and } w \in W\}$, and let $H$ be the subgroup generated by $T_X$. 
Then $T \cap H = T_X$.}

\bigskip\noindent
Let $\Omega$ be the (regular) graph defined as follows. The set of vertices of $\Omega$ is $S$. Two 
vertices $s,t \in S$ are joined by an edge in $\Omega$ if $m_{s\,t} < +\infty$ and $m_{s\,t}$ is odd. 
Let $\Omega_1, \dots, \Omega_m$ be the connected components of $\Omega$ and, for $1 \le i\le m$, let 
$Y_i$ be the set of vertices of $\Omega_i$. We suppose that $X \cap Y_i \neq \emptyset$ for $1 \le i\le 
p$ and $X \cap Y_i = \emptyset$ for $p+1 \le i\le m$. Let $C_2= \{ \pm 1\}$ be the cyclic group of 
order 2. There is an epimorphism $\kappa: W \to C_2$ which sends $s$ to 1 for all $s \in Y_1 \cup \dots 
\cup Y_p$ and $s$ to $-1$ for all $s \in Y_{p+1} \cup \dots \cup Y_m$. We have $\kappa(t)=1$ for 
all $t \in T_X$, thus $\kappa (H)=\{1\}$, and $\kappa (t)=-1$ for all $t \in T \setminus T_X$. This 
shows that $T \cap H=T_X$.

\bigskip\noindent
{\bf Assertion 4.} {\it Let $X \subset S$, $X \neq \emptyset$, such that $W_X$ is finite. Then $N_W(W_X)$ has 
infinite index in $W$.}

\bigskip\noindent
Suppose that $N_W(W_X)$ has finite index in $W$. Since $W_X$ is finite, this implies that $Z_W(W_X)$ has 
finite index in $W$, too. Set $T_X= \{ wxw^{-1}; x \in X \text{ and } w \in W\}$. 
Then the fact that $Z_W(W_X)$ has finite index in $W$ implies that $T_X$ is finite.
Let $H$ be the subgroup of $W$ generated by $T_X$. By 
definition, $H$ is normal. The subgroup $H$ is also a reflection group thus, by Proposition 4.2, there 
exists a subset $Y \subset H \cap T$ such that $(H,Y)$ is a Coxeter system. We have $H \cap T = 
T_X$ (Assertion 3) and $T_H= \{hyh^{-1}; y \in Y \text{ and } h \in H\} \subset T_X$, thus $T_H$ is 
finite. By Proposition 4.2, it follows that $H$ is finite. As pointed out before, $H$ is also normal, 
and this contradicts Assertion 2.

\bigskip\noindent
{\bf Assertion 5.} {\it Let $H$ be a non-trivial finite subgroup of $W$. Then $N_W(H)$ has infinite 
index in $W$.}

\bigskip\noindent
Let $\Pc (H)$ denote the intersection of all the parabolic subgroups that contain $H$. As before, $\Pc 
(H)$ is a non-trivial finite parabolic subgroup. By Assertion 4, it follows that $N_W(\Pc (H))$ has 
infinite index in $W$. Since $N_W(H) \subset N_W( \Pc (H))$, we conclude that $N_W(H)$ has infinite 
index in $W$.
\qed 

\bigskip\noindent
{\bf Proof of Theorem 4.1.} We start with the second part of the theorem. We assume that $\Gamma$ is 
connected, non-spherical, and non-affine. Let $G$ be a finite index subgroup of $W$ 
and let $A,B$ be two subgroups of $G$ such that $G = A \times B$. Take an essential element $w \in W$ 
(a Coxeter element, for example). Since $G$ has finite index in $W$, there is some non-trivial power $w^m$ of $w$ 
which belongs to $G$. By Corollary 2.5, the element $w^m$ is also essential, thus we can assume that $w \in 
G$. We write $w= w_A w_B$ where $w_A \in A$ and $w_B \in B$. The element $w$ has infinite order thus 
either $w_A$ or $w_B$ has infinite order (say $w_A$ has infinite order). The subgroup $\langle w \rangle$ has 
finite index in $Z_W(w)$ (by Theorem 2.6), $w_A \in Z_W(w)$, and $w_A$ has infinite order, thus there 
exist $p,q \in \Z \setminus \{ 0\}$ such that $w^p=w_A^q$. By Corollary 2.5, it follows that $w_A$ is 
essential. The group $\langle w_A \rangle$ has finite index in $Z_W(w_A)$ (by Theorem 2.6) and $\langle 
w_A \rangle \times B \subset Z_W(w_A)$, thus $B$ is finite. Finally, from the inclusion $G \subset N_W(B)$ and 
the fact that $G$ has finite index in $W$, follows that $N_W(B)$ has finite index in $W$, hence, by 
Proposition 4.3, $B=\{1\}$.

\bigskip\noindent
We turn now to the first part of Theorem 4.1. By the above, it suffices to consider the case where 
$\Gamma$ is connected and of affine type. Let $A,B$ be two subgroups of $W$ such that $W=A \times B$. By 
\cite{Bourb1}, Chap. 6, \S\ 2, there exist a finite Coxeter group $W_0$ and an irreducible 
representation $\rho: W_0 \to GL(\Z^{n-1})$ such that $W= \Z^{n-1} \rtimes W_0$. (Irreducible means in 
this context that the corresponding linear representation $W_0 \to GL(\R^{n-1})$ is irreducible. Note also 
that $n$ is the rank of $(W,S)$, but this fact is not needed for our purpose.) One of the components, $A$ or 
$B$, must be infinite (say $A$). The subgroup $\Z^{n-1}$ has finite index in $W$, thus $A \cap \Z^{n-
1}$ has finite index in $A$, therefore $A \cap \Z^{n-1} \neq \{1\}$ since $A$ is infinite. The subgroup 
$A \cap \Z^{n-1}$ is normal since both, $A$ and $\Z^{n-1}$ are normal, thus $(A \cap \Z^{n-1}) \otimes 
\R$ is a non-trivial linear subspace of $\R^{n-1}$ invariant by the action of $W_0$, therefore $(A \cap 
\Z^{n-1}) \otimes \R= \R^{n-1}$ since $\rho$ is irreducible. This implies that $A \cap \Z^{n-1}$ has 
finite index in $\Z^{n-1}$ and, consequently, in $W$.  
Thus $A$ has finite index in $W$. We conclude
that $B$ is finite and hence, by 
Proposition 4.3, $B=\{1\}$.
\qed


\section{Remak decompositions}

We start this section with the statement of the Krull-Remak-Schmidt theorem, and refer to
\cite{Robin1}, Section 3.3, for the proof.

\bigskip\noindent
{\bf Theorem 5.1 (Krull-Remak-Schmidt).} {\it 
\begin{enumerate}
\item
Let $G$ be a group which satisfies either the maximal or the minimal condition on normal subgroups.
Then $G$ has a Remak decomposition.
\item
Let $G$ be a group which satisfies both, the maximal and the minimal conditions on normal subgroups.
If $G=G_1 \times \dots \times 
G_p$ and $G=H_1 \times \dots \times H_q$ are two Remak decompositions, then $p=q$ and there exist a 
permutation $\sigma \in \Sym_p$ and a central automorphism $\varphi: G \to G$ such that $\varphi(G_i)= 
H_{\sigma (i)}$ for all $1 \le i\le p$. 
\end{enumerate}}

\bigskip\noindent
We turn now to the main result of the section. Let $\Gamma$ be a Coxeter graph and let $(W,S)$ be the 
Coxeter system of type $\Gamma$. Write $W=W_{X_1} \times \dots \times W_{X_b} \times W_{Z_3}$, where 
$W_{X_1}, \dots, W_{X_b}$ are the non-spherical irreducible components of $W$ and $W_{Z_3}$ is the finite part. 
By Theorem 5.1, we can choose a Remak decomposition $W_{Z_3}=H_1 \times \dots \times H_q$ of $W_{Z_3}$ which is 
unique up to a central automorphism. Then, by Theorem 4.1, $W=W_{X_1} \times \dots \times W_{X_b} 
\times H_1 \times \dots \times H_q$ is a Remak decomposition.

\bigskip\noindent
{\bf Theorem 5.2.} {\it Let $W=G_1 \times \dots \times G_m$ be a Remak decomposition of $W$. Then 
$m=b+q$ and there exist a permutation $\sigma \in \Sym_m$ and a central automorphism $\varphi: W \to W$ 
such that $\varphi(W_{X_i})= G_{\sigma(i)}$ for all $1 \le i \le b$ and $\varphi(H_j)= G_{\sigma( b+j)}$ 
for $1 \le j\le q$.}

\bigskip\noindent
{\bf Corollary 5.3.} {\it The group $\Aut(W_{X_1}) \times \dots \times \Aut( W_{X_b}) \times \Aut 
(W_{Z_3})$ has finite index in $\Aut(W)$.}

\bigskip\noindent
In many cases, for example when $\Gamma$ has no infinite labels (see \cite{HoRoTa1}), the group 
$\Out(W)$ is finite. However, there are Coxeter groups having infinite outer automorphism groups, for 
example the free product of $n$ copies ($n\ge 3$) of $C_2= \{ \pm 1\}$.

\bigskip\noindent
{\bf Corollary 5.4.} {\it Let $(\tilde W, \tilde S)$ be a Coxeter system. Write $\tilde W= \tilde 
W_{\tilde X_1} \times \dots \times \tilde W_{\tilde X_d} \times \tilde W_{\tilde {Z_3}}$, where $\tilde 
W_{\tilde X_1}, \dots, \tilde W_{\tilde X_d}$ are the non-spherical irreducible components of $\tilde W$ and 
$\tilde W_{\tilde Z_3}$ is the finite part. Then $W$ and $\tilde W$ are isomorphic if and only if 
$W_{Z_3}$ and $\tilde W_{\tilde Z_3}$ are isomorphic, $b=d$, and there exists a permutation $\sigma 
\in \Sym_b$ such that $W_{X_i}$ is isomorphic to $\tilde W_{\tilde X_{\sigma (i)}}$ for all $1 \le 
i\le b$.}

\bigskip\noindent
Call a Coxeter group {\it rigid} if it cannot be defined by two non-isomorphic Coxeter graphs. 

\bigskip\noindent
{\bf Corollary 5.5.} {\it The Coxeter group $W$ is rigid if and only if $W_{X_1}, \dots, W_{X_b}, 
W_{Z_3}$ are all rigid.}

\bigskip\noindent
Rigid Coxeter groups include, notably, those Coxeter groups that act effectively, properly, and 
cocompactly on contractible manifolds (see \cite{ChaDav1}). By a recent result of Caprace and 
M\"uhlherr \cite{CapMuh1}, the infinite Coxeter groups defined by connected Coxeter graphs that have no 
infinite labels are also rigid. Note that, by Corollary 5.5, this last result extends to the Coxeter 
groups defined by disjoint unions of non-spherical connected Coxeter graphs having no infinite labels. 
On the other hand, an interesting construction of non-rigid (and irreducible) Coxeter groups is given 
in \cite{BrMcMuNe1}. We refer to \cite{Muhlh1} for a recent exposition on rigidity in Coxeter groups 
and on the so-called isomorphism problem.

\bigskip\noindent
We start now the proof of Theorem 5.2 with the following lemma which will be also used in the next 
two sections.

\bigskip\noindent
For a group $G$ and a subset $E \subset G$, we denote by $Z_G(E)= \{w \in G; wx=xw \text{ for all } x 
\in E\}$ the centralizer of $E$ in $G$.

\bigskip\noindent
{\bf Lemma 5.6.} {\it Let $G$ be a group and let $A,B$ be two subgroups of $G$ such that $G=A \times 
B$. Let $E$ be a subset of $G$. Then $Z_G(E)= (Z_G(E) \cap A) \times (Z_G(E) \cap B)$.}

\bigskip\noindent
{\bf Proof.} Take $w=(w_A,w_B) \in Z_G(E)$. Let $x=(x_A,x_B) \in E$. We have $wx= (w_Ax_A, w_Bx_B)= xw= 
(x_Aw_A, x_Bw_B)$, thus $w_A$ and $x_A$ commute, therefore $w_Ax= (w_Ax_A, x_B)= (x_Aw_A, x_B)= xw_A$. 
This shows that $w_A \in Z_G(E)$ (and also that $w_B \in Z_G(E)$).
\qed

\bigskip\noindent
{\bf Proof of Theorem 5.2.} Recall that $W_{Z_3}= H_1 \times \dots \times H_q$ and $W= W_{X_1} \times 
\dots \times W_{X_b} \times H_1 \times \dots \times H_q$ are Remak decompositions. Recall also that the 
center of $W$ is included in $W_{Z_3}$. In particular, $Z(W)= Z(H_1) \times \dots \times Z(H_q)$. Now, 
we assume given a Remak decomposition $W=G_1 \times \dots \times G_m$.

\bigskip\noindent
Take $i \in \{1, \dots, b\}$ and set $\tilde W_i= W_{X_i} \times Z(W)= W_{X_i} \times Z(H_1) \times 
\dots \times Z(H_q)$. The group $\tilde W_i$ is the centralizer in $W$ of $W_{X_1} \times \dots \times 
W_{X_{i-1}} \times W_{X_{i+1}} \times \dots \times W_{X_b} \times W_{Z_3}$, thus, by Lemma 5.6,
\begin{equation}\label{eq1}
\tilde W_i= (\tilde W_i \cap G_1) \times \dots \times (\tilde W_i \cap G_m)\,.
\end{equation}
The inclusion $Z(W) \subset \tilde W_i$ implies that $Z(G_j) \subset (\tilde W_i \cap G_j)$ for all $1 
\le j \le m$. Hence, the quotient of the equality \eqref{eq1} by $Z(W)$ gives the isomorphism
\begin{equation}\label{eq2}
W_{X_i} \simeq (( \tilde W_i \cap G_1) / Z(G_1)) \times \dots \times (( \tilde W_i \cap G_m) / 
Z(G_m))\,.
\end{equation}
Since $W_{X_i}$ is indecomposable, it follows that there exists $\chi(i) \in \{ 1, \dots, m\}$ such that
\begin{equation}\label{eq3}
\begin{split}
(\tilde W_i \cap G_{\chi(i)})/ Z(G_{\chi(i)}) &\simeq W_{X_i}\,,\\
(\tilde W_i \cap G_j) /Z(G_j) &=\{1\} \quad \text{for } j \neq \chi(i)\,.
\end{split}\end{equation}
For $1 \le j\le m$, we denote by $\kappa_j: W \to G_j$ the projection on the $j$-th component. By 
\eqref{eq3}, the restriction $\kappa_{\chi(i)}: W_{X_i} \to G_{\chi (i)}$ of $\kappa_{\chi(i)}$ to $W_{X_i}$ 
is injective, and $\kappa_j(W_{X_i}) \subset Z(G_j)$ for $j \neq \chi(i)$. Set $\hat W_{X_i}= W_{X_1} \times 
\dots \times W_{X_{i-1}} \times W_{X_{i+1}} \times \dots \times W_{X_b} \times W_{Z_3}$ and define 
$\varphi_i: W= W_{X_i} \times \hat W_{X_i} \to W$ by
\[
\varphi_i(w)= \left\{
\begin{array}{ll}
w\quad&\text{if }w \in \hat W_{X_i}\,,\\
w \cdot {\displaystyle\prod_{j \neq \chi(i)} \kappa_j(w)^{-1}} \quad &\text{if } w \in W_{X_i}\,.
\end{array}\right.
\]
Then $\varphi_i$ is a well-defined central automorphism, it sends $W_{X_i}$ into $G_{\chi(i)}$, and it 
restricts to the identity on $\hat W_{X_i}$. Moreover, $\varphi_i$ and $\varphi_k$ commute if $i \neq 
k$. Set
\[
\varphi= \prod_{i=1}^b \varphi_i\,.
\]
Then $\varphi: W \to W$ is a central automorphism and $\varphi(W_{X_i}) \subset G_{\chi (i)}$ for all $1 
\le i\le b$. So, upon replacing $W_{X_i}$ by $\varphi(W_{X_i})$ for all $1 \le i\le b$ if necessary, we 
can assume that $W_{X_i} \subset G_{\chi(i)}$ for all $1 \le i\le b$.

\bigskip\noindent
For a group $G$ and two subgroups $A,B$ of $G$, we denote by $[A,B]$ the subgroup of $G$ generated by 
$\{ \alpha^{-1} \beta^{-1} \alpha \beta; \alpha \in A \text{ and } \beta \in B\}$. In particular, the 
equality $[A,B]=\{1\}$ means 
that every element of $A$ commutes with every element of $B$.

\bigskip\noindent
From the equality $W=W_{X_i} \times \hat W_{X_i}$ we deduce that $G_{\chi(i)}= \kappa_{\chi(i)} (W)= 
W_{X_i} \cdot \kappa_{\chi(i)} (\hat W_{X_i})$ and $[ W_{X_i}, \kappa_{\chi(i)} (\hat W_{X_i})] = 
\{1\}$. Moreover, $W_{X_i} \cap \kappa_{\chi(i)} (\hat W_{X_i})$ lies in the center of $W_{X_i}$ which 
is trivial, thus $W_{X_i} \cap \kappa_{\chi(i)} (\hat W_{X_i}) =\{1\}$, therefore $G_{\chi(i)} = 
W_{X_i} \times \kappa_{\chi(i)} (\hat W_{X_i})$. Since $G_{\chi(i)}$ is indecomposable, it follows that 
$\kappa_{\chi(i)} (\hat W_{X_i}) = \{1\}$ and $G_{\chi(i)} = W_{X_i}$. It also follows that $\chi(i) 
\neq \chi(k)$ if $i \neq k$. So, up to a permutation of the $G_j$'s, we can assume that $W_{X_i}=G_i$ 
for all $1 \le i\le b$.

\bigskip\noindent
At this point, we have the Remak decompositions $W=W_{X_1} \times \dots \times W_{X_b} \times H_1 
\times \dots \times H_q$ and $W=W_{X_1} \times \dots \times W_{X_b} \times G_{b+1} \times \dots \times 
G_m$. We also have $W_{Z_3}= H_1 \times \dots \times H_q$. Let $\pi: W \to W_{Z_3}= H_1 \times \dots 
\times H_q$ be the projection on $W_{Z_3}$. We have $\pi(G_{b+j}) 
\simeq G_{b+j}$ for all $1 \le j\le m-b$, and $W_{Z_3}= \pi(G_{b+1}) \times \dots \times \pi(G_m)$ is a 
Remak decomposition. By Theorem 5.1, it follows that $m=b+q$ and, up to a central automorphism and a 
permutation of the $G_{b+j}$'s, $\pi(G_{b+j})=H_j$ for all $1 \le j\le q$. The equality 
$\pi(G_{b+j})=H_j$ implies that $G_{b+j} \subset W_{X_1} \times \dots \times W_{X_b} \times H_j$, and 
the decomposition $W=W_{X_1} \times \dots \times W_{X_b} \times G_{b+1} \times \dots \times G_m$ 
implies that $G_{b+j} \subset Z_W(W_{X_1} \times \dots \times W_{X_b}) = H_1 \times \dots \times H_q$, 
thus $G_{b+j} \subset (W_{X_1} \times \dots \times W_{X_p} \times H_j) \cap (H_1 \times \dots \times 
H_q)= H_j$, that is, $G_{b+j}=H_j$ (since $\pi(G_{p+j})=H_j$).
\qed 


\section{Virtual strong Remak decompositions}

Let $\Gamma$ be a Coxeter graph and let $(W,S)$ be the Coxeter system of type $\Gamma$. Write $W= 
W_{X_1} \times \dots \times W_{X_a} \times W_{Z_2} \times W_{Z_3}$, where $W_{X_1}, \dots, W_{X_a}$ are 
the non-spherical and non-affine irreducible components of $W$, $W_{Z_2}$ is the affine part, and $W_{Z_3}$ is 
the finite part. We choose a finite index subgroup $R$ of $W_{Z_2}$ isomorphic to $\Z^d$, where $d$ is 
the affine dimension of $W$, and a decomposition $R=R_1 \times \dots \times R_d$ of $R$ as a direct 
product of $d$ copies of $\Z$. Then, by Theorem 4.1, $G=W_{X_1} \times \dots \times W_{X_a} \times R_1 
\times \dots \times R_d$ is a virtual strong Remak decomposition of $W$.

\bigskip\noindent
{\bf Theorem 6.1.} {\it Let $H=H_1 \times \dots \times H_m$ be a virtual strong Remak decomposition of 
$W$. Then $m=a+d$ and there exist
\begin{itemize}
\item
a virtual strong Remak decomposition $K=A_1 \times \dots \times A_a \times B_1 \times \dots \times B_d$ 
of $W$,
\item
a central automorphism $\varphi: K \to K$,
\item
a permutation $\sigma \in \Sym_m$,
\end{itemize}
such that
\begin{itemize}
\item
$A_i$ is a finite index subgroup of $W_{X_i}$ and $\varphi (A_i)$ is a finite index subgroup of 
$H_{\sigma (i)}$ for all $1 \le i\le a$,
\item
$B_j$ is isomorphic to $\Z$, $\varphi (B_j)$ is a finite index subgroup of $H_{\sigma (a+j)}$ for 
all $1 \le j\le d$, and $B=B_1 \times \dots \times B_d$ is a finite index subgroup of $R$.
\end{itemize}}

\bigskip\noindent
{\bf Corollary 6.2.} {\it Every Coxeter group has a unique virtual strong Remak decomposition up to equivalence.}

\bigskip\noindent
{\bf Corollary 6.3.} {\it Let $(\tilde W, \tilde S)$ be a Coxeter system. Write $\tilde W= \tilde 
W_{\tilde X_1} \times \dots \times W_{\tilde X_c} \times \tilde W_{\tilde Z_2} \times \tilde 
W_{\tilde Z_3}$, where $\tilde W_{\tilde X_1}, \dots, \tilde W_{\tilde X_c}$ are the non-spherical and 
non-affine irreducible components of $\tilde W$, $\tilde W_{\tilde Z_2}$ is the affine part, and $\tilde 
W_{\tilde Z_3}$ is the finite part. Then $W$ and $\tilde W$ are commensurable if and only if $a=c$, the 
affine dimension of $W$ is equal to the affine dimension of $\tilde W$, and there exists a permutation 
$\sigma \in \Sym_a$ such that $W_{X_i}$ and $\tilde W_{\tilde X_{\sigma(i)}}$ are commensurable for all 
$1 \le i\le a$.}

\bigskip\noindent
The proof of Theorem 5.2 uses the fact that the center of any infinite irreducible Coxeter group 
is trivial, but the proof of Theorem 6.1 needs the following stronger result.

\bigskip\noindent
{\bf Proposition 6.4.} {\it Assume $\Gamma$ to be connected, non-spherical, and non-affine. 
Let $G$ be a finite index subgroup of $W$. Then $Z_W(G)= \{1\}$.}

\bigskip\noindent
{\bf Proof.} Take $w_0 \in Z_W(G)$. Suppose first that $w_0$ has finite order. Let $H= \langle w_0 
\rangle$ be the subgroup of $W$ generated by $w_0$. The group $G$ is included in $N_W(H)$ and $G$ has 
finite index in $W$, thus $N_W(H)$ has finite index in $W$, therefore, by Proposition 4.3, $H=\{1\}$ 
and $w_0=1$.

\bigskip\noindent
Now, suppose that $w_0$ has infinite order. Take some essential element $w \in W$. Since $G$ has finite 
index in $W$, there is some non-trivial power $w^m$ of $w$ which belongs to $G$. So, upon replacing $w$ by $w^m$, 
we may assume that $w \in G$. The subgroup $\langle w \rangle$ has finite index in $Z_W(w)$ (by Theorem 
2.6), $w_0\in Z_W(w)$, and $w_0$ has infinite order, thus there exist $p,q \in \Z \setminus \{0\}$ such 
that $w_0^p=w^q$. By Corollary 2.5, this implies that $w_0$ is essential. Finally, $G \subset 
Z_W(w_0)$, $G$ has finite index in $W$, and $\langle w_0 \rangle$ has finite index in $Z_W(w_0)$, thus 
$\langle w_0 \rangle$ is a finite index subgroup of $W$. As pointed out in the introduction, such a 
Coxeter group, being virtually $\Z$, must be 
the Coxeter group of type $\tilde A_1$ which is of affine type: a contradiction.
\qed

\bigskip\noindent
{\bf Proof of Theorem 6.1.} First, observe the following fact: if $G$ is a group, $H$ is a finite index 
subgroup of $G$, and $K$ is any subgroup, then $H \cap K$ is a finite index subgroup of $K$.

\bigskip\noindent
Now, let $G=W_{X_1} \times \dots \times W_{X_a} \times R$, where $W_{X_1}, \dots, W_{X_a}$ are the 
non-spherical and non-affine irreducible components of $W$, and $R$ is a finite index subgroup of the affine 
part isomorphic to $\Z^d$. Let $H=H_1 \times \dots \times H_m$ be a virtual strong Remak decomposition.

\bigskip\noindent
Set $H_i'= H_i \cap G$ for all $1 \le i \le m$ and $H'= H_1' \times \dots \times H_m'$. Since $G$ has 
finite index in $W$, each $H_i'$ has finite index in $H_i$, thus $H'$ has finite index in $H$ (and in 
$W$). So, upon replacing $H$ by $H'$ and $H_i$ by $H_i'$ for all $1 \le i\le m$, we can assume that $H 
\subset G$.

\bigskip\noindent
Set $A_i= W_{X_i} \cap H$ for all $1 \le i \le a$, $B=R \cap H$, and $K=A_1 \times \dots \times A_a 
\times B$. The group $A_i$ is a finite index subgroup of $W_{X_i}$ for all $1 \le i\le a$ and $B$ is a 
finite index subgroup of $R \simeq \Z^d$ (in particular, $B \simeq \Z^d$), thus $K$ is a finite index 
subgroup of both, $G$ and $H$.

\bigskip\noindent
Now, we show that $Z(H)=B$. Let $w \in Z(H)$. Since $H \subset G$, the element $w$ can be written in 
the form $w_1 \cdots w_a u$ where $w_i \in W_{X_i}$ for all $1 \le i\le a$ and $u \in R$. We have $w_i 
\in Z_{W_{X_i}} (A_i)$ and $Z_{W_{X_i}} (A_i)$ is trivial by Proposition 6.2, thus $w_i=1$ for all $1 
\le i \le a$ and $w=u \in R \cap H = B$. This shows that $Z(H) \subset B$. The reverse inclusion $B 
\subset Z(H)$ is obvious.

\bigskip\noindent
Take $i \in \{ 1, \dots, a\}$ and set $\tilde A_i= A_i \times B$ and $\hat A_i= A_1 \times \dots 
\times A_{i-1} \times A_{i+1} \times \dots \times A_a \times B$. Let $G'= W_{X_1} \times \dots \times 
W_{X_a} \times B$. 
It is easily shown using the equality $Z(H)= B= R \cap H$ that $H \subset G'$.
By Proposition 6.2, we have $Z_{G'} (\hat A_i) = W_{X_i} \times B$, thus $Z_H(\hat 
A_i)= (W_{X_i} \times B) \cap H$. Furthermore, since $B \subset H$, we have $(W_{X_i} \times B) \cap H 
= (W_{X_i} \cap H) \times B = A_i \times B = \tilde A_i$, hence $Z_H(\hat A_i) = \tilde A_i$. By Lemma 
5.6, it follows that
\begin{equation}\label{eq4}
\tilde A_i = (\tilde A_i \cap H_1) \times \dots \times (\tilde A_i \cap H_m)\,.
\end{equation}
The inclusion $B=Z(H) \subset \tilde A_i$ implies that $Z(H_j) \subset (\tilde A_i \cap H_j)$ for all 
$1 \le j\le m$. Hence, the quotient of the equality \eqref{eq4} by $B$ gives the isomorphism
\begin{equation}\label{eq5}
A_i \simeq (( \tilde A_i \cap H_1)/ Z(H_1)) \times \dots \times (( \tilde A_i \cap H_m)/ Z(H_m))\,.
\end{equation}
Since $W_{X_i}$ is strongly indecomposable and $A_i$ is a finite index subgroup of $W_{X_i}$, it follows 
that there exists $\chi(i) \in \{1, \dots, m\}$ such that
\begin{equation}\label{eq6}
\begin{split}
(\tilde A_i \cap H_{\chi(i)}) / Z(H_{\chi(i)}) &\simeq A_i\,,\\
(\tilde A_i \cap H_j) / Z(H_j) &=\{1\} \quad \text{for } j \neq \chi(i)\,.
\end{split}
\end{equation}
For $1 \le j\le m$, we denote by $\kappa_j: H \to H_j$ the projection on the $j$-th component. By 
\eqref{eq6}, the restriction $\kappa_{\chi(i)} : A_i \to H_{\chi(i)}$ of $\kappa_{\chi(i)}$ to $A_i$ is 
injective, and $\kappa_j(A_i) \subset Z(H_j)$ for $j \neq \chi(i)$. Define $\varphi_i: K= A_i \times 
\hat A_i \to K$ by
\[
\varphi_i (w)= \left\{
\begin{array}{ll}
w\quad &\text{if } w \in \hat A_i\,,\\
w \cdot {\displaystyle \prod_{j \neq \chi(i)} \kappa_j (w)^{-1}} \quad&\text{if } w \in A_i\,.
\end{array}\right.
\]
Then $\varphi_i$ is well-defined since, for $j \neq \chi(i)$, $\kappa_j(w) \in Z(H_j) \subset Z(H) = B 
\subset K$, for all $w \in A_i$. It is a central automorphism, it sends $A_i$ into $H_{\chi(i)}$, and it 
restricts to the identity on $\hat A_i$. Moreover, $\varphi_i$ and $\varphi_k$ commute if $i \neq k$. 
Set
\[
\varphi= \prod_{i=1}^a \varphi_i\,.
\]
Then $\varphi: K \to K$ is a central automorphism and $\varphi(A_i) \subset H_{\chi(i)}$ for all $1 \le 
i\le a$. So, upon replacing $A_i$ by $\varphi(A_i)$ for all $1 \le i\le a$, we can assume that $A_i 
\subset H_{\chi(i)}$ for all $1 \le i \le a$.

\bigskip\noindent
From the equality $K= A_i \times \hat A_i$ follows that $\kappa_{\chi(i)} (K) = A_i \cdot 
\kappa_{\chi(i)} (\hat A_i)$ and $[A_i, \kappa_{\chi(i)} (\hat A_i) ] = \{1\}$. Moreover, $A_i \cap 
\kappa_{\chi(i)} (\hat A_i)$ lies in the center of $A_i$ which is trivial by Proposition 6.2, thus $A_i 
\cap \kappa_{\chi(i)} (\hat A_i) = \{1\}$, therefore $\kappa_{\chi(i)} (K)= A_i \times \kappa_{\chi(i)} 
(\hat A_i)$. The group $K$ is a finite index subgroup of $H$, thus $\kappa_{\chi(i)} (K)$ is a finite 
index subgroup of $H_{\chi(i)}$. Now, $H_{\chi(i)}$ is strongly indecomposable, thus $\kappa_{\chi(i)} 
(\hat A_i) = \{1\}$, $A_i=\kappa_{\chi(i)} (K)$, and $A_i$ is a finite index subgroup of $H_{\chi(i)}$. It also 
follows that $\chi(i) \neq \chi(k)$ if $i \neq k$. So, up to a permutation of the $H_j$'s, we can 
assume that $A_i$ is a finite index subgroup of $H_i$ for all $1 \le i \le a$.

\bigskip\noindent
Recall that $B= Z(H)= Z(H_1) \times \dots \times Z(H_m)$. Furthermore, for $1 \le i \le a$, we have 
$Z(H_i) = \kappa_i(B) \subset \kappa_i (\hat A_i)= \{1\}$. So, $B= Z(H) = Z(H_{a+1}) \times \dots \times 
Z(H_m)$. In particular, $B$ is a subgroup of $H_{a+1} \times \dots \times H_m$. Consider the inclusions
\[
K= A_1 \times \dots \times A_a \times B \subset A_1 \times \dots \times A_a \times H_{a+1} \times \dots 
\times H_m \subset H_1 \times \dots \times H_a \times H_{a+1} \times \dots \times H_m=H\,.
\]
Since $K$ has finite index in $H$, $K=A_1 \times \dots \times A_a \times B$ has finite index in $A_1 
\times \dots \times A_a \times H_{a+1} \times \dots \times H_m$, thus $B$ has finite index in $H_{a+1} 
\times \dots \times H_m$. For $1 \le j\le m-a$, we set $B_j= B \cap H_{a+j} = Z(H_{a+j})$. The group 
$B_j$ has finite index in $H_{a+j}$ and $H_{a+j}$ is strongly indecomposable, thus $B_j$ is non-trivial 
and indecomposable. On the other hand, $B_j$ is a subgroup of $B$ which is a finitely generated free 
abelian group, therefore $B_j$ is isomorphic to $\Z$. Finally, $B=B_1\times \dots \times B_{m-a} \simeq \Z^d$, thus 
$m-a=d$.
\qed


\section{Finite irreducible Coxeter groups}

Recall that
the connected spherical type Coxeter graphs are precisely the graphs $A_n$ ($n\ge 1$), $B_n$ ($n \ge 
2$), $D_n$ ($n \ge 4$), $E_6$, $E_7$, $E_8$, $F_4$, $H_3$, $H_4$, $I_2(p)$ ($p \ge 5$). Here we 
use the notation $I_2(6)$ for the Coxeter graph $G_2$. We may also use the notation $I_2(3)$ for 
$A_2$, and $I_2(4)$ for $B_2$. We number the vertices of each of these graphs as it is usually 
done. The numbering coincides with the standard numbering of the simple roots if $W$ is a Weil 
group (see \cite{Humph2}, p. 58), and the vertices of $H_3$ and $H_4$ are numbered so that $s_1$ 
and $s_2$ are joined by an edge labeled by $5$.

\bigskip\noindent
Let $\Gamma$ be a spherical type connected Coxeter graph and let $(W,S)$ be the Coxeter system of 
type $\Gamma$. For $w \in W$, we denote by $\lg (w)$ the word length of $w$ with respect to $S$. 
The group $W$ has a unique element of maximal length, $w_0$, which satisfies $w_0^2=1$ and 
$w_0Sw_0=S$. The following proposition is a list of well-known properties on $W$ which can be 
found for instance in \cite{Bourb1}.

\bigskip\noindent
{\bf Proposition 7.1.}
{\it 
\begin{enumerate}
\item There exists a permutation $\theta : S \to S$ such that $w_0sw_0= \theta(s)$ for all $s \in S$, 
and $\theta^2= \Id$.
\item The group $N_W(S)= \{ w\in W; wSw^{-1}=S\}$ is cyclic of order 2 generated by $w_0$.
\item If $\theta= \Id$, then $Z(W)=N_W(S)$. If $\theta \neq 
\Id$, then $Z(W)= \{1\}$.
\item We have $\theta \neq \Id$ if and only if $\Gamma \in \{ A_n; n\ge 2\} \cup \{ D_n; n\ge 5 
\text{ and } n \text{ odd}\} \cup \{ I_2(p); p\ge 5 \text{ and } p \text{ odd}\} \cup \{E_6\}$.
\end{enumerate}}

\bigskip\noindent
{\bf Theorem 7.2.} {\it Let $\Gamma$ be a connected spherical type Coxeter graph and let $(W,S)$ be the 
Coxeter system of type $\Gamma$. Then $W$ is decomposable if and only if $\Gamma 
\in \{I_2(p); p\ge 6 \text{ and }p \equiv 2(\text{mod }4)\} \cup \{B_n; n\ge 3; \text{ and } n 
\text{ odd}\} \cup \{H_3,E_7\}$. 
In that case, a Remak decomposition of $W$ is isomorphic to $Z(W) \times W/Z(W)$.}

\bigskip\noindent
{\bf Remark.} If either $\Gamma= I_2(p)$ ($p\ge 6$ and $p\equiv 2(\text{mod }4)$) or $\Gamma= 
B_n$ ($n \ge 3$ and $n$ odd), then the factor $W/Z(W)$ is a Coxeter group. It is the 
Coxeter group of type $I_2({p \over 2})$ if $\Gamma= I_2(p)$ ($p \ge 6$ and $p \equiv 2( 
\text{mod }4)$) and it is the Coxeter group of type $D_n$ if $\Gamma=B_n$ ($n \ge 3$ and $n$ 
odd). If $\Gamma \in \{ H_3, E_7\}$, then $W/Z(W)$ is not a Coxeter group. It is the alternating 
group on 5 letters if $\Gamma=H_3$ and it is $SO_7(2)$ if $\Gamma=E_7$.

\bigskip\noindent
{\bf Proof of Theorem 7.2.} We shall often use centralizers of parabolic subgroups (of finite 
Coxeter groups) in our proof. These can be easily computed with the techniques of \cite{Howle1}, 
\cite{Deodh1}, \cite{BriHow1}, or \cite{Borch1}, together with a suitable computer program. For instance, 
we used the package ``Chevie'' of GAP.

\bigskip\noindent
Let $\Gamma$ be a connected finite type Coxeter graph and let $(W,S)$ be the Coxeter system of 
type $\Gamma$. We assume given two subgroups $A,B \subset W$ such that $W=A \times B$, and we 
argue case by case.

\bigskip\noindent
{\bf Case $\Gamma= I_2(p)$, $p \ge 3$.} Left to the reader.

\bigskip\noindent
{\bf Case $\Gamma= A_n$, $n \ge 3$.} The case $n=3$ is left to the reader. So, we assume $n \ge 
4$. Let $X= \{s_3,s_4, \dots, s_n\}$. We have $Z_W(X)= \{1,s_1\}$. By Lemma 5.6, if follows that 
either $s_1 \in A$ or $s_1 \in B$ (say $s_1 \in A$). We conclude that $A=W$ since $A$ is normal 
and $W$ is normally generated by $s_1$.

\bigskip\noindent
{\bf Case $\Gamma= B_n$, $n \ge 3$.} Let $C_2= \{ \pm 1\}$ be the cyclic group of order $2$. Then 
$W= (C_2)^n \rtimes \Sym_n$, where $\Sym_n$ denotes the $n$-th symmetric group. Let $w_0$ be the 
longest element of $W$. Then $w_0 = [-1,-1, \dots, -1] \in (C_2)^n$, and $Z(W)= \{1,w_0\}$.

\bigskip\noindent
Consider the natural projection $\eta: W \to \Sym_n$. We have $\Sym_n= \eta(A) \cdot \eta(B)$, 
and $[\eta(A), \eta(B)]= \{1\}$. Notice that $\eta(A) \cap \eta(B) \subset Z(\Sym_n)= \{1\}$, thus 
we actually have $\Sym_n= \eta(A) \times \eta(B)$. By the previous case, it follows that either $\eta(A)= 
\{1\}$ or $\eta(B)=\{1\}$ (say $\eta(B)= \{1\}$).

\bigskip\noindent
Assume $B \neq \{1\}$. Let $\gamma: (C_2)^n \to C_2$ be the epimorphism defined by $\gamma 
[\varepsilon_1, \dots, \varepsilon_n] = \prod_{i=1}^n \varepsilon_i$, and let $K= \Ker \gamma$. 
There are precisely four normal subgroups of $W$ contained in $(C_2)^n$: $\{1\}$, $Z(W)=\{1, 
w_0\}$, $K$, and $(C_2)^n$. Thus either $B=Z(W)$, or $B=K$, or $B=(C_2)^n$. Observe that 
$\eta(A)= \Sym_n$, thus there exists $u \in (C_2)^n$ such that $u \cdot (1,2) \in A$. We cannot 
have either $B=K$ or $B=(C_2)^n$, because $v=[1,-1,-1,1, \dots, 1] \in K \subset (C_2)^n$, and 
$v$ does not commute with $u \cdot (1,2)$. So, $B=Z(W)= \{1,w_0\}$.

\bigskip\noindent
It remains to show that $n$ is odd. The quotient of $W$ by $A$ determines an epimorphism $\mu: W 
\to C_2$ which satisfies $\mu(w_0)=-1$, and such an epimorphism exists only if $n$ is odd.

\bigskip\noindent
{\bf Case $\Gamma= D_n$, $n \ge 4$.} Let $K$ be the subgroup of $(C_2)^n$ defined above. We have 
$W=K \rtimes \Sym_n$. Then one can easily prove, using the same arguments as in the case 
$\Gamma=B_n$, that either $A=\{1\}$ or $B=\{1\}$. (Note that $[-1,-1, \dots, -1] \not \in K$ if $n$ 
is odd.)

\bigskip\noindent
{\bf Case $\Gamma= H_3$.} Let $w_0$ be the longest element of $W$. Then $Z(W) = \{1, w_0\}$. Let 
$\gamma: W \to \{ \pm 1\}$ be the epimorphism defined by $\gamma(s_1) = \gamma (s_2) = \gamma 
(s_3) =-1$, and let $K= \Ker \gamma$. Then $\gamma(w_0)=-1$ and $W=K \times Z(W)$. Note also that 
$W$ is normally generated by $s_1$, and $K$ is normally generated by $s_1w_0$.

\bigskip\noindent
The group $Z_W (s_1)$ is generated by $\{ s_1, s_3, w_0\}$, and isomorphic to $(C_2)^3$. By Lemma 5.6, 
we have $Z_W(s_1)= (Z_W(s_1) \cap A) \times (Z_W(s_1) \cap B)$. Observe that $|Z_W(s_1)|=8$, 
thus $|Z_W(s_1) \cap (A \cup B)| \ge 5$, therefore $\{ s_1, s_1w_0, s_3, sè_3w_0\} \cap (A \cup 
B) \neq \emptyset$. Furthermore, $s_1$ and $s_3$ are conjugate in $W$, thus there exists $a \in 
\{0,1\}$ such that either $s_1w_0^a \in A$ or $s_1w_0^a \in B$ (say $s_1w_0^a \in A$). We 
conclude that $K \subset A$, thus either $A=K$ or $A=W$.

\bigskip\noindent
{\bf Case $\Gamma=H_4$.} Let $X= \{ s_3,s_4\}$. Set
\[
v= s_2s_1s_3 s_2s_1s_2 s_1s_4s_3 s_2s_1s_2 s_1s_3s_2 s_1s_2s_3 s_4s_3s_2\,.
\]
Then $Z_W(X)$ is generated by $\{s_1,v\}$ and has the presentation $\langle s_1,v | s_1^2= v^2 = 
(s_1v)^6=1 \rangle \simeq W_{I_2(6)}$. Let $w_0$ be the longest element of $W$. Then $Z(W)=\{1, 
w_0\}$ and $(s_1v)^3=(vs_1)^3=w_0$. By Lemma 5.6, we have $Z_W(X)= (Z_W(X) \cap A) \times (Z_W(X) 
\cap B)$. By the case $\Gamma= I_2(6)$ treated before, it follows that there exists $a \in 
\{0,1\}$ such that either $s_1w_0^a \in A$ or $s_1w_0^a \in B$ (say $s_1w_0^a \in A$). We 
conclude that $A=W$ since $W$ is normally generated by $s_1w_0^a$ (whatever is $a \in \{ 0,1 
\}$).

\bigskip\noindent
{\bf Case $\Gamma=F_4$.} Let $X=\{s_1,s_2\}$, let $Y= \{ s_1,s_2,s_3\}$, and let $w_Y$ be the 
longest element of $W_Y$. The group $Z_W(X)$ is generated by $\{ s_4, w_Y\}$, and has the 
presentation $\langle s_4, w_Y| s_4^2= w_Y^2= (s_4w_Y)^6=1 \rangle \simeq W_{I_2(6)}$. Let $w_0$ 
be the longest element of $W$. Then $Z(W)= \{1,w_0\}$ and $(s_4w_Y)^3=(w_Ys_4)^3=w_0$. By Lemma 
5.6, we have $Z_W(X)= (Z_W(X) \cap A) \times (Z_W(X) \cap B)$. By the case $\Gamma=I_2(6)$ 
treated before, it follows that there exists $a \in \{0,1\}$ such that either $s_4w_0^a \in A$ or 
$s_4 w_0^a \in B$ (say $s_4 w_0^a \in A$). The generator $s_3$ is conjugate to $s_4$, thus we 
also have $s_3w_0^a \in A$. Similarly, there exists $b \in \{0,1\}$ such that either $s_1w_0^b, 
s_2 w_0^b \in A$, or $s_1w_0^b, s_2w_0^b \in B$. We cannot have $s_1w_0^b, s_2 w_0^b \in B$, 
since $s_2w_0^b$ and $s_3w_0^a$ do not commute, thus $s_1w_0^b, s_2w_0^b \in A$. We conclude that 
$w_0= (s_1s_2s_3s_4)^6 \in A$, thus $s_1,s_2,s_3,s_4 \in A$, therefore $A=W$.

\bigskip\noindent
{\bf Case $\Gamma=E_6$.} The centralizer of $X=\{ s_1,s_2,s_3,s_4\}$ in $W$ is $\{1,s_6\}$, thus, 
by Lemma 5.6, either $s_6 \in A$, or $s_6 \in B$ (say $s_6 \in A$). We conclude that $A=W$ since 
$W$ is normally generated by $s_6$.

\bigskip\noindent
{\bf Case $\Gamma=E_7$.} Let $w_0$ be the longest element of $W$. Then $Z(W)= \{1, w_0\}$. Let 
$\gamma: W \to C_2= \{ \pm 1\}$ be the epimorphism defined by $\gamma(s_i)=-1$ for all $1 \le 
i\le 7$, and let $K= \Ker \gamma$. Then $\gamma(w_0)=-1$ and $W=K \times Z(W)$. Note also that 
$W$ is normally generated by $s_1$, and $K$ is normally generated by $s_1w_0$.

\bigskip\noindent
Using the same argument as in the case $\Gamma=H_4$, but with the centralizer of $X= \{s_2, 
s_3,s_5, s_6, s_7\}$, we show that there exists $a \in \{0,1\}$ such that either $s_1w_0^a \in A$ 
or $s_1w_0^a \in B$ (say $s_1w_0^a \in A$). This implies that $K \subset A$, thus either $A=K$ or 
$A=W$.

\bigskip\noindent
{\bf Case $\Gamma=E_8$.} This case can be handle in the same way as the case $\Gamma=H_4$, but 
with the centralizer of $X=\{s_2,s_4,s_5,s_6,s_7,s_8\}$.
\qed



\bigskip\bigskip\noindent
{\bf Luis Paris},

\smallskip\noindent 
Institut de Math\'ematiques de Bourgogne, UMR 5584 du CNRS, Universit\'e de Bourgogne, B.P. 
47870, 21078 Dijon cedex, France

\smallskip\noindent
E-mail: {\tt lparis@u-bourgogne.fr}


\end{document}